\documentclass[12pt]{article}
\usepackage{amsthm,amsfonts,amssymb,amscd}

\begin{document}


\begin{center}
\large \bf Birational geometry of singular Fano double spaces\\
of index two
\end{center}\vspace{0.5cm}

\centerline{A.V.Pukhlikov}\vspace{0.5cm}

\parshape=1
3cm 10cm \noindent {\small \quad\quad\quad \quad\quad\quad\quad
\quad\quad\quad {\bf }\newline In this paper we describe the birational geometry of Fano double spaces
$V\stackrel{\sigma}{\to}{\mathbb P}^{M+1}$ of index 2 and dimension
$\geqslant 8$ with at mostquadratic singularities of rank
$\geqslant 8$, satisfying certain additional conditions of general position: we prove that these varieties have no structures of a rationally connected fibre space over a base of dimension $\geqslant 2$, that every birational map $\chi\colon V\dashrightarrow V'$ onto the total space of a Mori fibre space
$V'/{\mathbb P}^1$ induces an isomorphism $V^+\cong V'$ of the blow up
$V^+$ of the variety $V$ along $\sigma^{-1}(P)$, where $P\subset
{\mathbb P}^{M+1}$ is a linear subspace of codimension 2, and that every birational map of the variety $V$ onto a Fano variety $V'$ with ${\mathbb
Q}$-factorial terminal singularities and Picard number 1 is an isomorphism. We give an explicit lower estimate for the codimension of the set of varieties $V$ with worse singularities or not satisfying the conditions of general position, quadratic in $M$. The proof makes use of the method of maximal singularities and the improved $4n^2$-inequality for the self-intersection of a mobile linear system.\vspace{0.1cm}

Bibliography: 20 items.} \vspace{1cm}

14E05, 14E07\vspace{0.1cm}

\vspace{0.3cm}

{\bf 1. Statement of the main result.} The aim of the present paper is to extend the results of \cite{Pukh10} to singular double spaces of index two. At the same time, we generalize and improve these results and make them more precise. Let ${\mathbb P}={\mathbb P}^{M+1}$ be the complex projective space, where $M\geqslant 7$. Hypersurfaces of degree $2M$ in ${\mathbb P}$ are parameterized by points of the projective space
$$
{\cal W}={\mathbb P}(H^0({\mathbb P},{\cal O}_{\mathbb P}(2M))).
$$
If the singular set of a hypersurface $W\subset {\mathbb
P}$ of degree $2M$ is of codimension $\geqslant 4$ in ${\mathbb
P}$, then the double cover $\sigma\colon V\to {\mathbb P}$, branched over $W$, is an irreducible factorial variety. If singularities of $V$ are terminal, then $V$ is a Fano variety of index two:
$$
\mathop{\rm Pic}V=\mathop{\rm Cl}V={\mathbb Z}H,\quad K_V=-2H,
$$
where $H=\sigma^*H_{\mathbb P}$ is the $\sigma$-pull back of the class of a hyperplane in ${\mathbb P}$. Such varieties are realized as hypersurfaces
$$
y^2=F(x_0,\dots,x_{M+1})
$$
in the weighted projective space ${\mathbb P}(1,\dots,1,M)={\mathbb P}(1^{M+2},M)$, where the equation $f(x_*)=0$ is the equation of the hypersurface $W$, and $y$ is the coordinate of weight $M$. If $P\subset{\mathbb P}$ is a linear subspace of codimension 2, then the projection $\alpha_P\colon{\mathbb P}\dashrightarrow{\mathbb P}^1$ from that subspace induces on $V$ the structure
$$
\pi_P=\alpha_P\circ\sigma\colon V\dashrightarrow{\mathbb P}^1
$$
of a Fano-Mori fibre space, the fibres of which are Fano double spaces of index 1. Let us consider the integer-valued function
$$
\beta\colon{\mathbb Z_{\geqslant 7}=\{M\,|\,M\geqslant 7\}}\to{\mathbb Z},
$$
defined in the following way: $\beta(7)=3$ and for $M\geqslant 8$
$$
\beta(M)=\frac12(M-5)(M-4)+1.
$$
A rationally connected fibre space means, as usual, a surjective morphism $\lambda\colon Y\to S$ of projective varieties, such that a fibre of general position $\lambda^{-1}(s),s\in S$ and the the base $S$ are rationally connected.

Here is the main result of the present paper.

{\bf Theorem 1.} {\it There is a Zariski open subset ${\cal U}\subset {\cal W}$, such that
$$
\mathop{\rm codim}(({\cal W}\backslash {\cal U})\subset {\cal W})\geqslant\beta(M)
$$
and the double cover $\sigma\colon V\to{\mathbb P}$, branched over any hypersurface $W\in {\cal U}$, is an irreducible reduced factorial variety with at most terminal singularities and satisfying the following two properties:

{\rm (i)} for any rationally connected fibre space $\lambda\colon Y\to S$ over a positive-dimensional base $S$ and any birational map
$\chi\colon V\dashrightarrow Y$ onto the total space $Y$ (if there are such maps) the equality $S={\mathbb P^1}$ holds and for some isomorphiam $\beta\colon{\mathbb P}^1\to S$ and some linear subspace $P\subset{\mathbb P}$ of cdimension 2 the diagram
$$
\begin{array}{rcccl}
& V & \stackrel{\chi}{\dashrightarrow} & Y &\\
\phantom{M}\pi_P & \downarrow & & \downarrow & \lambda\phantom{M}\\
& {\mathbb P}^1 & \stackrel{\beta}{\to} & S &\\
\end{array}
$$
commutes, that is, $\lambda\circ\chi=\beta\circ\pi_P$;

{\rm (ii)} every birational map $\chi\colon V\dashrightarrow V'$ onto a Fano variety $V'$ with ${\mathbb Q}$-factorial terminal singularities and Picard number 1 is a biregular isomorphism.}

Theorem 1 immediately implies the standard set of facts about birational geometry of double covers $V\to{\mathbb P}$, branched over a hypersurface $W\in {\cal U}$.

{\bf Corollary 1.} {\it For every variety $V$ which is the double cover of the projective space ${\mathbb P}$, branched over some hypersurface $W\in {\cal U}$, the following claims are true.

{\rm (i)} On the variety $V$ there are no structures of a rationally connected fibre space (and, therefore, of a Fano-Mori fibre space) over a base of dimension $\geqslant 2$. In particular, on $V$ there are no structures of a conic bundle and del Pezzo fibration, and the variety $V$ itself is non-rational.\vspace{0.1cm}

{\rm (ii)} Assume that there is a birational map
$\chi\colon V\dashrightarrow V'$ onto the total space of a Mori fibre space $\pi'\colon V'\to S'$ with $\mathop{\rm dim} S'\geqslant 1$. Then $S'={\mathbb P}^1$ and for some subspace $P\subset {\mathbb P}$ of codimension 2 the birational map $\chi^{-1}\colon V'\dashrightarrow V$ is the blow up of the subvariety $\sigma^{-1}(P)$ (in particular, $\chi^{-1}$ is regular), and moreover, for some isomorphism $\beta\colon
{\mathbb P}^1\to S$ the equality
$$
\pi_P\circ\chi^{-1}=\beta^{-1}\circ\pi'
$$
holds.

{\rm (iii)} The groups of biregular and birational automorphisms of the variety $V$ are equal: $\mathop{\rm Bir}V=\mathop{\rm Aut}V$.
}

{\bf Proof of the corollary.} The first claim follows from the part (i) of Theorem 1 in an obvious way. Similarly, the claim (iii) follows immediately from the part (ii) of Theorem 1.

Let us consider the claim (ii). By Theorem 1, (i) we have $S'={\mathbb P}^1$ and for some linear subspace $P\subset{\mathbb P}$ of codimension 2 the birational map $\chi$ transforms the fibres of the projection $\pi_P$ into the fibres of the projection $\pi'$. The fact that $\chi$ induces an isomorphism of the blow up $V_P$ of the variety $V$ along $\sigma^{-1}(P)$ and the variety $V'$ is shown in Section 3. Q.E.D.\vspace{0.3cm}


{\bf 2. Conditions of general position.}  The open subset ${\cal U}\subset {\cal W}$ consists of hypersurfaces $W\subset{\mathbb P}$, satisfying the conditions that will be now defined. Let
$$
F(x_0,x_1,\dots,x_{M+1})=0
$$
be the equation of the hypersurface $W$ in ${\mathbb P}={\mathbb P}^{M+1}$.

(R0) For any subspace $P\subset{\mathbb P}$ of codimension 3 the restriction $F|_P$ is not a square of a polynomial (of degree $M$).

Now let $p\in W$ be a point and $z_1,\dots,z_{M+1}$ a coordinate system on an affine subset ${\mathbb A}^{M+1}\subset{\mathbb P}$, where $p=(0,\dots,0)$. Let us write down the affine equation of the hypersurface $W$ with respect to this system of coordinates:
$$
f=q_1+q_2+\dots+q_{2M},
$$
where $q_i(z_*)$ are homogeneous of degree $i\geqslant 1$. The point $p$ is non-singular on $W$ if and only if $q_1\not\equiv 0$.

(R1) Assume that $p\in W$ is a non-singular point. The the inequality
$$
\mathop{\rm rk}q_2|_{\{q_1=0\}}\geqslant 4
$$
holds.

(R2) Assume that $p\in W$ is a singular point. The the inequality
$$
\mathop{\rm rk}q_2\geqslant 7
$$
holds.

The set ${\cal U}\subset{\cal W}$ consists of hypersurfaces, satisfying the global condition (R0), the local condition (R1) at every non-singular point and the local condition (R2) at every singular point.

Let $W\in{\cal U}$ be an arbitrary hypersurface. Consider the double cover $\sigma\colon V\to{\mathbb P}$, branched over $W$. By the condition (R2), the inequality
\begin{equation} \label{25.09.2019.1}
\mathop{\rm codim\,}(\mathop{\rm Sing}W\subset{\mathbb P})\geqslant 7
\end{equation}
holds, so that
$$
\mathop{\rm codim\,}(\mathop{\rm Sing}V\subset V)\geqslant 7
$$
and for that reason $V$ is an irreducible reduced factorial (by Grothendieck's theorem \cite{CL}) variety. In \cite[ï. 2.1]{Pukh15a} it was shown that the condition to have at most quadratic singularities of rank $\geqslant r$ is stable under blow ups, so that the singularities of the variety $V$, which has at most quadratic singularities of rank $\geqslant 8$, are terminal.

{\bf Theorem 2.} {\it The inequality
$$
\mathop{\rm codim\,}(({\cal W}\backslash{\cal U})\subset{\cal W})\geqslant\beta(M)
$$
holds.}

{\bf Proof:} elementary computations. It is easy to check that the codimension of the closed subset ${\cal W}_0\subset{\cal W}$ of hypersurfaces that do not satisfy the condition (R0) is
$$
{3M-2\choose M-2}-{2M-2\choose M-2}-3(M-1),
$$
the codimension of the closed subset ${\cal W}_1\subset{\cal W}$ of hypersurfaces violating the condition (R1) at at least one non-singular point is at least
$$
\frac12M(M-7)+3,
$$
and the codimension of the closed subset ${\cal W}_2\subset{\cal W}$ of hypersurfaces violating the condition (R2) at at least one singular point is at least
$$
\frac12(M-5)(M-4)+1.
$$
It is easy to check that the function $\beta(M)$, defined above, gives the minimum of these three expressions for $M\geqslant 8$. Q.E.D. for Theorem 2.

Theorem 2 implies that Theorem 1 follows from the claim below.

{\bf Theorem 3.} {\it The double cover $\sigma\colon V\to{\mathbb P}$, branched over a hypersurface $W\subset{\cal U}$, satisfies the claims (i) and  (ii) of Theorem 1.}

The rest of this paper contains a proof of Theorem 3.\vspace{0.3cm}


{\bf 3. The method of maximal singularities.} Both claims (i) and (ii) of Theorem 1 and the part of claim (ii) of Corollary 1 which is yet to be shown, are proved by means of the method of maximal singularities \cite[Chapter 2]{Pukh13a}. Starting from this moment, we fix the variety $V$, which is the double cover of the projective space ${\mathbb P}$, branched over a hypersurface $W\in{\cal U}$.

Let $\lambda\colon Y\to S$ be a rationally connected fibre space over a positive-dimensional base and $\Sigma_Y$ a linear system of divisors on $Y$, which is the $\lambda$-pull back of some very ample linear system on the base $S$. Furthermore, let $V'$ be a Fano variety with ${\mathbb Q}$-factorial terminal singularities and Picard number 1 and $\Sigma'=|-m'K_V'|$ a very ample linear system, where $m'\gg 1$ is large enough. Assume that there is a birational map
$$
\chi\colon V\dashrightarrow X,
$$
where either $X=Y$, or $X=V'$.

If $X=Y$, set $\Sigma$ to be the strict transform of $\Sigma_Y$ on $V$ with respect to $\chi$. If $X=V'$, set $\Sigma$ to be the strict transform of
$\Sigma'$ with respect to $\chi$. In any case, replacing, if necessary,
$\Sigma_Y$ or $\Sigma'$ by their symmetric square, we may assume that
$$
\Sigma\subset|-nK_V|=|2nH|
$$
is a mobile plurianticanonical linear system on the variety $V$. Recall that a prime exceptional divisor $E^*$ on some blow up $\widetilde{V}\stackrel{\mu}{\to}V$ of the variety $V$ is called a {\it maximal singularity} of the linear system $\Sigma$, if the {\it Noether-Fano inequality} holds:
$$
\mathop{\rm ord}\nolimits_{E^*}\mu^*\Sigma>n\cdot a(E^*),
$$
where $a(E^*)=a(E^*,V)$ is the discrepancy of $E^*$ with respect to $V$. The following fact is well known.

{\bf Theorem 4.} (i) {\it If $X=Y$ is the total space of a rationally connected fibre space over a positive-dimensional base $S$, then the linear system $\Sigma$ has a maximal singularity.}

(ii) {\it If $X=V'$ is a Fano variety with Picard number 1 and $n>m'$, then the linear system $\Sigma$ has a maximal singularity.}

For a {\bf proof}, see \cite[Chapter 2]{Pukh13a}.

The existence of a maximal singularity for the linear system $\Sigma$ can be expressed in the following way: the pair $(V,\frac{1}{n}\Sigma)$ is not canonical.

The following claim is of key importance in the proof of Theorem 3.

{\bf Theorem 5.} {\it Assume that the linear system $\Sigma\subset|2nH|$ has a maximal singularity. Then there exists a linear subspace $P\subset{\mathbb P}$ of codimension 2, satisfying the inequality}
\begin{equation}\label{19.09.2019.1}
\mathop{\rm mult}\nolimits_{\sigma^{-1}(P)}\Sigma>n.
\end{equation}

{\bf Proof} of Theorem 5 is given in Sections 5-8. Assuming that the theorem is shown, let us consider the birational map $\chi\colon V\dashrightarrow X$, where, if $X=V'$ is a Fano variety, we assume that, moreover, $n>m'$. Let $P\subset{\mathbb P}$ be the linear subspace of codimension 2, satisfying the inequality (\ref{19.09.2019.1}). By the symbol $|H_{\mathbb P}-P|$ we denote the pencil of hyperplanes in ${\mathbb P}$, containing $P$, and by the symbol $|H-P|$ the subsystem of the linear system $|H|$, consisting of divisors, containing $\sigma^{-1}(P)$. Let $\varphi\colon V^+\to V$ be the blow up of the subvariety $\sigma^{-1}(P)$. The symbol $E_P$ stands for the exceptional divisor of the blow up $\varphi$. Obviously, $\varphi$ resolves the singularities of the pencil $|H-P|$ and $\pi=\pi_P\circ\varphi\colon V^+\to{\mathbb P}^1$ is a morphism, the fibres of which are isomorphic to the corresponding divisors in the pencil $|H-P|$.

{\bf Proposition 1.} {\it The variety $V^+$ and every fibre of the projection $\pi$ have at most quadratic singularities of rank $\geqslant 5$. In particular, these varieties are factorial and have terminal singularities.}

{\bf Proof} is obtained by simple computations: in the system of affine coordinates $z_1,\dots,z_{M+1}$ on ${\mathbb A}^{M+1}\subset{\mathbb P}$ with the origin at the point $p\in P\cap W$ the variety $V$ is locally given by the equation
$$
y^2=q_1(z_*)+q_2(z_*)+\dots,
$$
if the point $p\in W$ is non-singular, and by the equation
$$
y^2=q_2(z_*)+\dots,
$$
if the point $p\in W$ is a singularity. If $R\in|H_{\mathbb P}-P|$ is an arbitrary hyperplane in the pencil, then $\sigma^{-1}(R)\in|H-P|$ is given, respectively, by the equation
$$
y^2=q_1(z_*)|_R+q_2(z_*)|_R+\dots
$$
or
$$
y^2=q_2(z_*)|_R+\dots\, .
$$
If the point $p\in W$ is non-singular, then $\sigma^{-1}(R)$ has a singularity at the point $o=\sigma^{-1}(p)$ in one case only: when $R=T_pW$. By the condition (R1), in that case the point $o=\sigma^{-1}(R)$ is a quadratic singularity of rank $\geqslant 5$.

If the point $p\in W$ is singular, then by the condition (R2) the rank of the quadratic form $q_2(z_*)|_R$ is at least 5, so that the point $o=\sigma^{-1}(R)$ is a quadratic singularity of rank $\geqslant 6$.

Therefore, all fibres of the projection $\pi\colon V^+\to{\mathbb P}$ have at most quadratic singularities of rank $\geqslant 5$, as we claimed. This implies, that the variety $V^+$ has this property, too (which can be also checked directly by means of the explicit formulas of for the blow up $\varphi$).

The second claim of Proposition 1 follows from the first one. Q.E.D. for the proposition.

The pull back of the divisorial class $\varphi^*H$ onto $V^+$ is denoted by the same symbol $H$. The fibre $\pi^{-1}(t)$ over a point $t\in{\mathbb P}^1$ is written as $F_t$, and its class in the Picard group as $F$. The canonical class of the variety $V^+$ we denote by the symbol $K^+$.

{\bf Proposition 2.} {\it The following equalities are true:
$$
\mathop{\rm Pic}V^+={\mathbb Z}H\oplus{\mathbb Z}E_P={\mathbb Z}K^+\oplus{\mathbb Z}F
$$
and} $K^+=-2H+E_P, F=H-E_P$.

{\bf Proof.} This is obvious.

Now let us consider the strict transform $\Sigma^+$ of the linear system $\Sigma$ on $V^+$. This is a mobile linear system, and moreover, for some $m\in{\mathbb Z_+}$ and $l\in{\mathbb Z}$ we have the inclusion
$$
\Sigma^+\subset|-mK^++lF|.
$$
Using Proposition 2, we compute $m$ and $l$:
$$
m=2n-\mathop{\rm mult}\nolimits_{\sigma^{-1}(P)}\Sigma,\quad l=2(\mathop{\rm mult}\nolimits_{\sigma^{-1}(P)}\Sigma-n)\geqslant 2.
$$
If $m=0$, then the linear system $\Sigma^+\subset|lF|$ is composed from the pencil $|F|$, which is equivalent to the claim (i) of Theorem 1. Thus, inorder to prove that claim in full, we assume that $m\geqslant 1$, and show that this assumption leads to a contradiction. In order to do that, we apply the method of maximal singularities to the variety $V^+$.

{\bf Proposition 3.} {\it The linear system $\Sigma^+$ has a maximal singularity: for some exceptional prime divisor $E^+$ over
$V^+$ the Noether-Fano inequality
$$
\mathop{\rm ord}\nolimits_{E^+}\Sigma^+>m\cdot a(E^+,V^+)
$$
holds, that is, the pair $(V^+,\frac{1}{m}\Sigma^+)$ is not canonical.}

{\bf Proof} is well known.

Now the proof of the claim (i) of Theorem 1 is completed almost in one line: let $F^*\in|F|$ be a fibre of the projection $\pi$, intersecting the centre of the maximal singularity $E^+$. The linear system $\Sigma^*=\Sigma^+|_{F^*}$ is non-empty, although it may no longer be mobile. Let $D^*\in\Sigma^*$ be a general divisor, $D^*\sim -mK_{F^*}$. The pair
$$
\left(F^*,\frac{1}{m}D^*\right)
$$
is not canonical. However, in \cite[Subsection 2.2]{Pukh05} it was shown that this is impossible (see also the proof of Theorem 1.2 in \cite{Pukh15a}). Here one has to take into account that $F^*$ is the double cover of the projective space ${\mathbb P}^M$, branched over the hypersurface $W^*\subset{\mathbb P}^M$ of degree $2M$, which is a hyperplane section of the hypersurface $W\subset{\mathbb P}$. If $p\in W^*$ is a point and $z_1,\dots,z_M$ a system of coordinates on ${\mathbb A}^M\subset{\mathbb P}^M$ with the origin at $p$, then the equation of the hypersurface $W^*$ can be written in the form
$$
f^*=q_1^*+q_2^*+\dots+q_{2M}^*,
$$
where $q^*_i(z_*)$ are homogeneous of degree $i\geqslant 1$. If the point $p\in W^*$ is non-singular, then the condition (R1) tells us that
$$
\mathop{\rm rk}q^*_2|_{\{q^*_1=0\}}\geqslant 2
$$
(the rank of a quadratic form can drop by at most 2 when it is restricted to a hyperplane). If the point $p\in W^*$ is singular, then the condition (R2) tells us that
$$
\mathop{\rm rk}q^*_2\geqslant 5
$$
(for the same reason). Therefore, the conditions of general position, imposed on the double space of index 1 in \cite{Pukh05}, are satisfied (they are precisely the conditions (W1) and (W) in \cite[\S 3]{Pukh15a}), and for that reason the global canonical threshold of the fibre $F^*$ is at least 1. Q.E.D. for the part (i) of Theorem 1.\vspace{0.3cm}


{\bf 4. Biregular refinements.} Now let us consider the claims (ii) of Theorem 1 and Corollary 1, which make make the ``roguh'' birational claim (i) of Theorem 1 stronger and more precise. These biregular facts are shown in \cite{Pukh2018b}. We do not reproduce the proof in full, just explain its main steps; the full details can be found in \cite{Pukh2018b}. Note, first of all, that if, in the notations of Theorem 4, the inequality $n\leqslant m'$ holds, then very well known arguments (see, for instance, the proof of \cite[Chapter 2, Proposition 1.6]{Pukh13a}) show that the map $\chi\colon V\dashrightarrow V'$ is a biregular isomorphism. So we may assume that $n>m'$, so that the system $\Sigma$ has a maximal singularity, and by Theorem 5 the inequality (\ref{19.09.2019.1}) holds. Let us consider again the linear system $\Sigma\subset|-mK^++lF|$. It can not be composed from a pencil, therefore $m\geqslant 1$. If $m>m'$, then, taking into account that $l\in{\mathbb Z}_+$, and arguing as, for instance, in \cite[Chapter 2, Proposition 1.2]{Pukh13a}, we obtain the claim of Proposition 3 and come to a contradiction in precisely the same way as in the relative case $X=Y$ (the end of Section 3). Therefore, we may assume that $m\leqslant m'$. In that case we argue in word for word the same way as in \cite[Subsection 1.4]{Pukh2018b}. This argument essentially repeats the proof of \cite[Chapter 2, Proposition 1.6]{Pukh13a} with the only difference: the map
$$
\chi^+=\chi\circ\varphi\colon V^+\dashrightarrow V'
$$
can not be an isomorphism, because $V^+$ and $V'$ have different Picard numbers,
$$
\rho(V^+)=2\neq 1=\rho(V').
$$
This completes the proof of Theorem 1.

Let us show the biregular refinement (ii) in Corollary 1. For the details, see \cite[Subsection 1.5]{Pukh2018b}, here we just describe the main steps of the arguments. So let $\pi'\colon V'\to S'$ be a Mori fibre space over a positive-dimensional base and $\chi\colon V\dashrightarrow V'$ a birational map. By Theorem 1, we have $S'={\mathbb P}^1$ and for some linear subspace $P\subset{\mathbb P}$ of codimension 2 the composition
$$
\chi^+=\chi\circ\varphi\colon V^+\dashrightarrow V'
$$
transforms the fibres of $\pi$ into the fibres of $\pi'$. Furthermore,
$$
\mathop{\rm Pic} V'\otimes{\mathbb Q}={\mathbb Q}K'\oplus{\mathbb Q}F'
$$
is a two-dimensional linear space over ${\mathbb Q}$, where $K'=K_{V'}$, and
$F'$ is the class of a fiber of the projection $\pi'$. Now we choose the linear system $\Sigma'$ on $V'$ in a different way: not pulling it back from the base $S'$. Let
$$
\Sigma'=|-m'K'+l'F'|
$$
be a very ample {\it complete} linear system. Its strict transform on $V^+$
$$
\Sigma^+\subset|-mK^++lF|
$$
satisfies the condition $l\in{\mathbb Z}_+$, if $l'\in{\mathbb Z}_+$ is large enough (recall that the strict transform of the pencil of fibres $|F'|$ on
$V^+$ is the pencil of fibres $|F|$). Now, if $m>m'$, then the linear system $\Sigma^+$ has a maximal singularity, and we get a contradiction as above, restricting a general divisor of the system $\Sigma^+$ onto the fibre $F^*$, intersecting the centre of the maximal singularity. Therefore,
$m\leqslant m'$ and, arguing in word for word the same way as in \cite[Subsection 1.5]{Pukh2018b}, we show that $\chi^+$ is an isomorphism in codimension 1. This implies that
$$
\Sigma=|-mK^++lF|
$$
is a very ample {\it complete} linear system. For that reason $\chi^+\colon V^+\to V'$  is a biregular isomorphism. This completes the proof of the claim (ii) of Corollary 1.

The remaining part of the paper is a proof of the key Theorem 5.\vspace{0.3cm}


{\bf 5. Maximal subvarieties of codimension 2.} Fix a mobile linear system $\Sigma\subset|2nH|$ with a maximal singularity $E^*$ on some blow up $\mu\colon\widetilde{V}\to V$ of the double space $V$. The method of proving Theorem 5 is the standard one: assume that the inequality (\ref{19.09.2019.1}) is {\it not} satisfied for all linear subspaces $P$ of codimension 2, that is,
$$
\mathop{\rm mult}\nolimits_{\sigma^{-1}(P)}\Sigma\leqslant n.
$$
Our aim is to bring this assumption to a contradiction. Let $B=\mu(E^*)$ be the centre of the maximal singularity on $V$.

{\bf Proposition 4.} {\it The following inequality holds:}
$$
\mathop{\rm codim}(B\subset V)\geqslant 3.
$$

{\bf Proof.} Assume the converse:
$$
\mathop{\rm codim}(B\subset V) = 2.
$$
We know that $\sigma(B)\subset{\mathbb P}$ is not a linear subspace of codimension 2, because
$$
\mathop{\rm mult}\nolimits_B\Sigma> n.
$$
Let us check that the arguments that were used in \cite{Pukh10} for excluding maximal singularities of that type work for singular double spaces, too. Somewhat abusing the notations, we use the symbol $P$ for a 6-dimensional linear subspace of of general position in ${\mathbb P}$. Since $P\cap\mathop{\rm Sing} W = \emptyset$, we have:
$$
V(P)=\sigma^{-1}(P)
$$
is a non-singular 6-dimensional variety, which is a double cover of $P\cong{\mathbb P}^6$, branched over a non-singular hypersurface $W(P)=W\cap P$. For the numerical Chow group of cycles of codimension 2 we have
$$
A^2V(P)={\mathbb Z}H^2_P,
$$
where $H_P$ is the $\sigma$-pull back of the class of a hyperplane in $P$, so that also
$$
A^2V={\mathbb Z}H^2,
$$
and for that reason $B\sim H^2$ or $2H^2$ or $3H^2$. Set $\overline{B}=\sigma(B)\subset{\mathbb P}$. By assumption, $\mathop{\rm deg}\overline{B}\geqslant 2$.

If $B=\sigma^{-1}(\overline{B})$, then the proof of Proposition 2.1 in \cite[Subsection 2.1]{Pukh10} works as it is: a general secant line $L$ of the subvariety $B$ does not intersect the set $\mathop{\rm Sing}W$. Therefore,
$$
\sigma^{-1}(\overline{B})=B\cup B',
$$
where $\mathop{\rm deg}\overline{B}\in\{2,4,6\}$ and $B\neq B'$. The case
$\mathop{\rm deg}\overline{B}=2$, when $\overline{B}$ is an irreducible quadric in the hyperplane $\langle\overline{B}\rangle\subset{\mathbb P}$, is excluded by word for word the same arguments as in \cite[Subsection 2.2]{Pukh10}, because a general 2-plane in $\langle\overline{B}\rangle$ does not intersect the set $\mathop{\rm Sing} W$.

Assume that $\mathop{\rm deg}\overline{B}=6$. Here the arguments of \cite[Subsection 2.3]{Pukh10} work with minimal modifications: let $C\subset\overline{B}$ be an irreducible curve, which is not contained in
$W$, and moreover $C\cap\mathop{\rm Sing}W=\emptyset$. We get the inequality
$$
\mathop{\rm mult}\nolimits_{B'}\Sigma > \frac{M-2}{M-1}\cdot n
$$
(in \cite{Pukh10} the dimension of the variety $V$ and the projective space
${\mathbb P}$ is equal to the value of the parameter $M$, in this paper $\mathop{\rm dim} V = M+1$). Looking at the self-intersection $Z=(D_1\circ D_2)$ of the linear system $\Sigma$, where $D_1,D_2\in\Sigma$ are general divisors, we get
$$
8n^2=\mathop{\rm deg}Z\geqslant 6((\mathop{\rm mult}\nolimits_B\Sigma)^2+
(\mathop{\rm mult}\nolimits_{B'}\Sigma)^2)>
6\cdot\left(1+\frac{(M-2)^2}{(M-1)^2}\right)n^2.
$$
Already for $M=5$ this inequality is impossible. Recall that in this paper
$M\geqslant 7$. The case $\mathop{\rm deg}\overline{B}=6$ is excluded.

Assume that $\mathop{\rm deg}\overline{B}=4$. Here we argue as in \cite[Subsection 2.4]{Pukh10}: Proposition 2.4 of that paper remains true and the proof of Proposition 2.5 works word for word in the cases 1), 2) and 3). Let us consider the last case 4), where the proof needs some minor modifications. In \cite[Subsection 2.4]{Pukh10} the symbol $P$ means a general linear subspace of dimension 3 in ${\mathbb P}$, so that $B_P=\overline{B}\cap P$ is a curve in ${\mathbb P}^3$. This is a curve of degree 4 that has in the case 4) a unique double point, so that the subvariety $\overline{B}$ contains a linear subspace $\Pi\subset{\mathbb P}$ of double points, $\mathop{\rm codim}(\Pi\subset{\mathbb P})=3$. We modify the proof given in \cite{Pukh10} in the following way.

Let us consider the net of hyperplanes in ${\mathbb P}$, containing $\Pi$. We denote this linear system by the symbol $|H_{\mathbb P}-\Pi|$. Let $\Theta\in|H_{\mathbb P}-\Pi|$ be a general element. Obviously,
$$
(\overline{B}\circ\Theta)=2\Pi+Q(\Theta),
$$
where $Q(\Theta)$ is an irreducible quadric in the linear subspace $\langle Q(\Theta)\rangle$ of codimension 2. This quadric is not contained in
$W$, because $\overline{B}\not\subset W$ by the Lefschetz theorem. Furthermore, the singular set of the hypersurface $W\cap\langle Q(\Theta)\rangle$ is of codimension at least 3 with respect to the subspace $\langle Q(\Theta)\rangle$, so that a 2-plane $\Lambda\subset\langle Q(\Theta)\rangle$ of general position does not intersect this set and the curve $W\cap\Lambda$ is non-singular. For that reason the proof of Lemma 2.1 in \cite{Pukh10} works. That lemma states that the surface $S=\sigma^{-1}(\Lambda)$ is contained in the base set of the linear system $\Sigma$. However, the planes $\Lambda$ sweep out a divisor on ${\mathbb P}$, so that the surfaces $S$ sweep out a divisor on $V$. This contradiction excludes the case $\mathop{\rm deg}B=4$ and completes the proof of Proposition 4.\vspace{0.3cm}


{\bf 6. Maximal singularities with the centre of codimension 3.} Assume now that the centre $B=\mu(E^*)$ of the maximal singularity on $V$ is of codimension 3. We need to bring this assumption to a contradiction. Set $\overline{B}=\sigma(B)$: this is a subvariety of codimension 3 in the projective space ${\mathbb P}$.

{\bf Proposition 5.} {\it The following inequality holds:} $\mathop{\rm deg} B \geqslant 2$.

{\bf Proof.} We have to exclude the option $\mathop{\rm deg}B=1$. If this equality were true, the subvariety $\overline{B}\subset{\mathbb P}$ would have been a linear subspace, and moreover,

(1) either $\sigma^{-1}(\overline{B})=B\cup B'$, where $B'=B$,

(2) or $\overline{B}\subset W$.

If the case (1) takes place, then the restriction $F|_{\overline{B}}$ of the equation of the hypersurface $W$ onto $\overline{B}$ is a full square, which is impossible by the condition (R0). If the case (2) takes place, then it is easy to check that
$$
\mathop{\rm codim}((\mathop{\rm Sing} W\cap\overline{B})\subset\overline{B})=3,
$$
so that
$$
\mathop{\rm codim}((\mathop{\rm
Sing} W\cap\overline{B})\subset{\mathbb P})=6,
$$
which contradicts the inequality (\ref{25.09.2019.1}). Q.E.D. for Proposition 5.

Now the proof of Proposition 3.1 in \cite{Pukh10} works word for word, excluding the option $\mathop{\rm codim}(B\cap V)=3$.\vspace{0.3cm}


{\bf 7. Maximal singularities with the centre of codimension $\geqslant 4$: the non-singular case.} Assume now that the centre $B=\mu(E^*)$ of the maximal singularity is of codimension $\geqslant 4$, and moreover,
$$
B\not\subset\mathop{\rm Sing}V,
$$
so that $\overline{B}=\sigma(B)\not\subset\mathop{\rm Sing}W$. This option is excluded by the arguments of \cite[\S 6]{Pukh10}, however, the fact that the hypersurface $W$ may be singular requires some fragments of the arguments to be modified. In the case $\overline{B}\not\subset W$ no changes are needed: the arguments of \cite[Subsection 6.1]{Pukh10} work word for word. So we assume that $\overline{B}\subset W$. Keeping the notations of \cite[\S 6]{Pukh10}, let $o\in B$ be a point of general position (in particular, $o\not\in\mathop{\rm Sing}V$) and $p=\sigma(o)\in\overline{B}$. Set
$$
\varphi\colon V^+\to V\quad\mbox{and}\quad\varphi_{\mathbb P}\colon{\mathbb P}^+\to{\mathbb P}
$$
to be the blow ups of the point $o$ and $p$, respectively (these notations were already used in Sections 3-4, but in the old sense they are no longer needed; we use these notations again because they compatible with the notations of \cite[\S 6]{Pukh10}). Let $E=\varphi^{-1}(o)$ and $E_{\mathbb P}=\varphi^{-1}_{\mathbb P}(p)$ be the exceptional divisors. The $8n^2$-inequality (see \cite[Proposition 4.1]{Pukh10} or \cite[Chapter 2, Section 4.1]{Pukh13a}) implies that there is a linear subspace $\Pi\subset E$ of codimension 2, satisfying the inequality
$$
\mathop{\rm mult}\nolimits_oZ+\mathop{\rm
mult}\nolimits_{\Pi}Z^+>8n^2
$$
(here $Z$, as above, is the self-intersection of the mobile linear system
$\Sigma$, and $Z^+$ its strict transform on $V^+$). By the symbol ${\mathbb T}$, as in \cite[Subsection 6.2]{Pukh10}, we denote the hyperplane ${\mathbb P}(T_PW)\subset E_{\mathbb P}$, which is naturally identified with the hyperplane ${\mathbb P}(T_o\sigma^{-1}(W))\subset E$. This identification will be assumed silently in the subsequent arguments. Recall that by our assumption $p\not\in\mathop{\rm Sing}W$ is a non-singular point on the branch divisor: it is for this reason that the proof given in \cite[\S 6]{Pukh10} works (with minor modifications).\vspace{0.3cm}

The morphism $\sigma$ induces a rational map of degree 2
$$
V^+\dashrightarrow{\mathbb P}^+,
$$
which {\it is not} a double cover. Its restriction onto the exceptional divisor
$$
\sigma_E\colon E\dashrightarrow E_{\mathbb P}
$$
is the projection $E\cong{\mathbb P}^M$ from a certain point $\xi\not\in{\mathbb E}\backslash{\mathbb T}$ onto the hyperplane ${\mathbb T}$, which is now considered as a hyperplane ${\mathbb T}\subset E_{\mathbb P}$ (see \cite[\S 6]{Pukh10}). Now the {\it simple case}, when $\xi\in\Pi$, so that $\Pi_{\mathbb P}=\sigma_E(\Pi)\subset{\mathbb T}$ is of codimension 2 in ${\mathbb T}$ and of codimension 3 in $E_{\mathbb P}$, is excluded by repeating word for word the arguments of \cite[Subsection 6.2]{Pukh10}. It remains to exclude the {\it hard case}, when $\xi\not\in\Pi$, so that $\Pi_{\mathbb P}=\sigma_E(\Pi)\subset{\mathbb T}$ is a hyperplane in ${\mathbb T}$ and a linear subspace of codimension 2 in $E_{\mathbb P}$. Here we argue as in \cite[Subsection 6.3]{Pukh10}. Let $\Lambda\subset{\mathbb P}$ be the unique subspace of codimension 2, such that $p\in\Lambda$ and $\Lambda^+\cap E_{\mathbb P}=\Pi_{\mathbb P}$, where $\Lambda^+\subset{\mathbb P}^+$ is the strict transform. By the assumption (R2) on the singularities of the hypersurface $W$ we have: $Q=\sigma^{-1}(\Lambda)$ is an irreducible subvariety.

{\bf Proposition 6.} {\it The strict transform $Q^+\subset V^+$ of the variety $Q$ does not contain $\Pi$.}

{\bf Proof.} Let $z_1,\dots,z_{M+1}$ be affine coordinates with the origin at the point $p$ and
$$
f=q_1(z_*)+q_2(z_*)+\dots+q_{2M}(z_*)
$$
the equation of the hypersurface $W$ in these coordinates. In a neighborhood of the point $o$ the local equation of the subvariety $Q\subset{\mathbb A}^M$ is of the form
$$
y^2=f|_{\Lambda}.
$$
Since $\Pi_{\mathbb P}\subset{\mathbb T}=\{q_1=0\}$, we get the equation
$$
y^2=q_2|_{\Lambda}+\dots+q_{2M}|_{\Lambda}
$$
where (denoting the affine part of the subspace $\Lambda$ by the same symbol) $\Lambda\subset\{q_1=0\}$. Therefore, $q_2|_{\Lambda}$ is the restriction of the quadratic form $q_2$ onto a hyperplane in the tangent hyperplane $\{q_1=0\}$. Therefore, the condition (R1) gives the inequality
$$
\mathop{\rm rk}q_2|_{\Lambda}\geqslant 2,
$$
so that the subvariety $Q$ has at the point $o$ a quadratic singularity of rank $\geqslant 3$. So the quadric $Q^+\cap E$ is a quadratic hypersurface of rank $\geqslant 3$ in its linear span $\langle Q^+\cap E\rangle$, which can not contain the linear subspace $\Pi$ (the subspace $\Pi$ is a hyperplane in $\langle Q^+\cap E\rangle$). Q.E.D. for the proposition.

Therefore, we get the equality
$$
\mathop{\rm mult}\nolimits_oQ+\mathop{\rm mult}\nolimits_{\Pi}Q^+=\mathop{\rm deg}Q=2,
$$
and the arguments of \cite[Subsection 6.3]{Pukh10} exclude the maximal singularity in the hard case.

We have excluded all options for the maximal singularity, if its centre $B$ is not contained in $\mathop{\rm Sing}V$.\vspace{0.3cm}


{\bf 8. Maximal singularities with the centre of codimension $\geqslant 4$: the singular case.} Assume now that $B\subset\mathop{\rm Sing}V$, so that $\overline{B}\subset\mathop{\rm Sing}W$, and moreover, $B$ has the maximal dimension among all centres of maximal singularities of the linear system
$\Sigma$. In particular, by the condition (R2) we have $\mathop{\rm codim}(B\subset V)\geqslant 7$.

This case is a serious challenge for the exclusion technique that was available in 2009 and used in \cite{Pukh10}. Let $Z=(D_1\circ D_2)$ be the self-intersection of the mobile system $\Sigma$, so that
$Z\sim 4n^2H^2$ and $\mathop{\rm deg}Z=8n^2$.

{\bf Proposition 7.} {\it The following inequality holds:}
$$
\mathop{\rm mult}\nolimits_BZ>8n^2.
$$

{\bf Proof.} Let $o\in B$ be a point of general position. It is sufficient to check that the germ of a quadratic singularity $o\in V$ satisifes the assumptions of the improved $4n^2$-inequality \cite[Section 2]{Pukh2017a}. Let $z_1,\dots,z_{M+1}$ be coordinates on the affine chart ${\mathbb A}^{M+1}\subset{\mathbb P}$, where $p=\sigma(o)=(0,\dots,0)$, so that the variety $V$ in a neighborhood of the point $o$ is the hypersurface
$$
\{y^2=q_2(z_*)+\dots+q_{2M}(z_*)\}\subset{\mathbb A}^{M+2}_{y,z_*}.
$$
Consider a general linear subspace (that is, an affine subspace, containing the point $o$) $P\subset{\mathbb A}^{M+2}$ of dimension 7. By the condition (R2) the point $o$ is an isolated quadratic singularity of the variety $V_P=V\cap P$ of maximal rank, so that for the blow up
$$
\varphi_P\colon V^+_P\to V_P
$$
of that point we have: $V^+_P$ is non-singular in a neighbourhood of the exceptional divisor $Q_P=\varphi_P^{-1}(o)$, which is a non-singular quadric in ${\mathbb P}^6$. Since the dimension $\mathop{\rm dim} B$ is maximal among the dimensions of all centres of maximal singularities of the system $\Sigma$, the subvariety $B$ is not strictly contained in another centre of a maximal singularity of this system. Therefore, the restriction $\Sigma_P=\Sigma|_{V_P}$ is a mobile linear system on $V_P$, and the pair $(V_P,\frac{1}{n}\Sigma_P)$ is canonical outside the point $o$ in a neighbourhood of that point, but has the point $o$ as the centre of a maximal singularity. For the self-intersection $Z_P$ of the system $\Sigma_P$ by   \cite[Section 1]{Pukh2017a} we have the inequality
$$
\mathop{\rm mult}\nolimits_oZ_P>4n^2\cdot 2=8n^2.
$$
However $Z_P=(Z\circ V_P)$ is the section of the cycle $Z$ by a general linear subspace $P\ni o$, so that
$$
\mathop{\rm mult}\nolimits_oZ>8n^2.
$$
This completes the proof of Proposition 7. Q.E.D.

In order to exclude the singular case, it remains to note that we got the inequality
$$
\mathop{\rm mult}\nolimits_BZ>\mathop{\rm deg}Z,
$$
which is impossible. This contradiction completes the proof of Theorem 5.\vspace{0.3cm}


{\bf 9. Historical remarks and acknowledgements.} The main theorem of \cite{Pukh10} was the first birational rigidity-type result for a large class of Fano varieties of index 2 and arbitrary $(\geqslant 5)$ dimension. In that paper the double spaces of index 2 were assumed to be non-singular. The paper \cite{Pukh10} was followed by a new breakthrough --- a description of birational geometry of general non-singular Fano hypersurfaces of index 2 and dimension $\geqslant 16$ in the projective space \cite{Pukh16a}. In \cite{Pukh2017a} the improved $4n^2$-inequality was shown for complete intersection singularities, which essentially simplified exclusion of maximal singularities, the centre of which is contained in the singular set of the variety. This made it possible to extend the results on the birational geometry of Fano varieties of index 2 to singular varieties and, in particular, to obtain an effective estimate for the codimension of the complement to the set of varieties covered by the main theorem. For Fano hypersurfaces of index 2 this work was done in \cite{Pukh2018b}, Fano double spaces of index 2 are studied in the present paper. Note that although the proof of the improved $4n^2$-inequality in its full generality makes use of the generalization of \cite[Proposition 5]{Pukh02b}, obtained in \cite{Suzuki15}, for varieties with hypersurface (in particular, quadratic) singularities the claim of \cite[Proposition 5]{Pukh02b} is sufficient.

The theorems of birational rigidity type, shown in the last ten years, as a rule, cover the families of Fano varieties and Mori fibre spaces with singularities of bounded type (see, for instance, \cite{Ahmad2014,Krylov2018,ChPark10,Foord2019,Pukh2019a} --- the list is far from being complete).

Double covers of the projective space play the role of a touchstone in this work from the very first steps of the theory of birational rigidity \cite{IM,I80,Pukh89a}. As the experience of \cite{EP2014,EP2018} shows, the quadratic singularities, the rank of which is bounded from below, form the most natural class of singularities for extending the  ``non-singular'' results to the ``singular'' context.

The author is grateful to the colleagues in the Divisions of Algebraic Geometry and Algebra of Steklov Institute of Mathematics for the interest to his work. The author is also grateful to the colleagues in the Algebraic Geometry research group at the University of Liverpool for the creative atmosphere and general support.

The author thanks The Leverhulme Trust for the support (Research
Project Grant RPG-2016-279).

\begin{flushleft}
Department of Mathematical Sciences,\\
The University of Liverpool
\end{flushleft}

\noindent{\it pukh@liverpool.ac.uk}

\end{document}